\documentclass[12pt]{article}
\usepackage[utf8]{inputenc}
\usepackage{amsmath}
\usepackage{amssymb}
\usepackage{amsfonts}
\def\End{\mathrm{End}}

\begin{document}

\begin{center} {\Large Rota-Baxter operators on cocommutative Hopf algebras}
\end{center}

\vspace{5mm}

\begin{center}
{\bf Maxim Goncharov}
\end{center}

\begin{abstract}
We generalize the notion of a Rota-Baxter operator on groups and the notion of a Rota-Baxter operator of weight 1 on Lie algebras and define and study the notion of a Rota-Baxter operator on a cocommutative Hopf algebra $H$.  If $H=F[G]$ is the group algebra of a group $G$ or $H=U(\mathfrak{g})$ the universal enveloping algebra  of a Lie algebra $\mathfrak{g}$, then we prove that Rota-Baxter operators on $H$ are in one to one correspondence with corresponding Rota-Baxter operators  on groups or Lie algebras.

\end{abstract}

\section{Introduction}

Given an arbitrary algebra $A$ over a field $F$ and a scalar $\lambda\in F$
a~linear operator $R\colon A\rightarrow A$ is called a Rota---Baxter operator
 on $A$ of weight~$\lambda$ if for all $x,y\in A$:
\begin{equation}\label{RB}
R(x)R(y) = R( R(x)y + xR(y) + \lambda xy )
\end{equation}
Then the pair $(A,R)$ is called a~Rota---Baxter algebra. If $R$ is a Rota-Baxter operator of weight $\lambda$ and $\alpha\in F$, then $\alpha R$ is a Rota-Baxter operator of weight $\alpha\lambda$. Thus, there are two principle cases: when $\lambda=0$ or $\lambda =1$.

Rota-Baxter operators for associative algebras first appear in the
paper of G. Baxter as a tool for studying integral operators  in the theory of probability and mathematical statistics
\cite{Br}. 

The combinatorial properties of (commutative) Rota-Baxter algebras
and operators were studied in papers of F.V. Atkinson, P. Cartier,
G.-C. Rota and the others (see \cite{Atk}-\cite{Car}). For basic
results and the main properties of Rota-Baxter algebras see
\cite{Guo}.

Independently, in early 80-th Rota-Baxter  operators on Lie algebras
naturally appear in papers of A.A. Belavin, V.G. Drinfeld \cite{BD}
and M.A. Semenov-Tyan-Shanskii \cite{STS} while studying the
solutions of the classical Yang-Baxter equation. It turns out that on quadratic Lie algebras skew-symmetric solutions of the classical Yang-Baxter equation are in one to one correspondence with skew-symmetric Rota-Baxter operators. 

If $\mathfrak{g}$ is a simple Lie algebra then  non-skew-symmetric $\mathfrak{g}$-invariant solutions of the classical Yang-Baxter equation (that sometimes called solutions of modified classical Yang-Baxter equation) on $\mathfrak{g}$ are in one to one correspondence with pairs $(R,B)$,  where $R$ is a Rota-Baxter operator of weight 1 satisfying $R+R^*+id =0$ and $B$ is a non-degenerate symmetric bilinear form on $\mathfrak{g}$ \cite{GME}. If $\mathfrak{g}$ is not simple, connections between non-skew-symmetric $\mathfrak{g}$-invariant solutions of the classical Yang-Baxter equation and Rota-Baxter operators were considered in \cite{GME1}. As a consequence of these results we can note, that every Lie biagrebra structure on a simple Lie algebra is induced by a Rota-Baxter operators of special type. 

When one considers the problem of quantization of a Lie bialgebra $(\mathfrak{g},\delta)$, one of the first step is to extend the comultiplication $\delta$ to a Poisson co-bracket on the universal enveloping algebra $U(\mathfrak{g})$ that can be done uniquely. From this point of view, it is natural to consider the question of extension of a Rota-Baxter operator $R$ from a Lie algebra $\mathfrak{g}$ to some reasonable operator on the universal enveloping algebra $U(\mathfrak{g})$. Unfortunately, it is not possible to extend $R$ to a Rota-Baxter operator of the algebra $U(\mathfrak{g})$ (that is, to a linear map  $B:U(\mathfrak{g})\mapsto U(\mathfrak{g})$ satisfying (\ref{RB})). 

Nevertheless, in \cite{OG} and \cite{FLMK} it was proved that a structure of a pre- or  a post-Lie algebra on a Lie algebra $\mathfrak{g}$ can be extended to some reasonable product on the universal enveloping algebra $U(\mathfrak{g})$. These results can be considered from the  point of view of Rota-Baxter operators: it turns out that a Rota-Baxter operator $R$ of weight $\lambda$ on  $\mathfrak{g}$ induces on $\mathfrak{g}$ a structure of a pre-Lie algebra (if $\lambda=0$) or a structure of a post-Lie algebra (if $\lambda\neq 0$). Here again we can ask if we can extend $R$ to an operator $B:U(\mathfrak{g})\mapsto U(\mathfrak{g})$ on the universal enveloping algebra $U(\mathfrak{g})$ in such a way that the extension of the pre-(or post-)Lie algebra structure on $U(\mathfrak{g})$ is somehow induced by $B$.

Recently, in \cite{GLY} it was introduced the notion of a Rota-Baxter operator (of weight 1) for groups. If $G$ is a group, then a map $B:G\mapsto G$ is called a Rota-Baxter operator on the group $G$ if for all $g,h\in G$:
$$
B(g)B(h)=B(gB(g)hB(g)^{-1}).
$$

A group $G$ with a Rota-Baxter operator $B$ is called a Rota-Baxter group.  In the same paper it was proved, that if $(G,B)$ is a Rota-Baxter Rota-Baxter Lie group, then the tangent map of $B$ at the identity is a Rota-Baxter operator of weight 1 on the Lie algebra of the Lie group $G$. Also, it was showed that many  results that are true for Rota-Baxter operators on algebras have corresponding analogs for Rota-Baxter operators on groups.

Lie algebras and groups can be regarded as foundations of two principle examples of cocommutative Hopf algebras. In this paper we in some sense combine notions of  Rota-Baxter operators of weight 1 on Lie algebras and of  Rota-Baxter operators on groups and give the definition of a Rota-Baxter operator (of weight 1) on cocommutative Hopf algebras. Note, that there already exist notions of Rota-Baxter of algebras and bialgebras (see \cite{JianZhang} and \cite{MaLiu}). These operators are different from the definition that we give. 

The paper organised as follows. In section 2 we give the definition of Rota-Baxter operator on a cocommutative Hopf algebra and  obtain some basic results  about it that are generalisations of known results for Rota-Baxter operators (of weight 1) on groups and algebras. In section 3, we consider two principle cases of cocommutative Hopf algebras - the universal enveloping algebra $U(\mathfrak{g})$ of a Lie algebra $\mathfrak{g}$ and the group algebra of a group $G$. We prove that Rota-Baxter operators on $U(\mathfrak{g})$ (resp. on $F[G]$) are in one-two-one correspondence with Rota-Baxter operators of weight 1 on $\mathfrak{g}$ (resp, on $G$). Given a Rota-Baxter operator $R$ of weight 1 on a Lie algebra $\mathfrak{g}$, one can define the structure of a post-Lie algebra on $\mathfrak{g}$. In section 4 we first show, that this extension  of a post-Lie algebra structure  to the universal enveloping algebra $U(\mathfrak{g})$ (that was found in \cite{FLMK}) can be defined using the Rota-Baxter Hopf operator $B:U(\mathfrak{g})\mapsto U(\mathfrak{g})$ that is the extension of $R$. Further, we prove that for a given arbitrary Rota-Baxter Hopf algebra $(H,B)$ one can define new multiplication $*$ and new antipod $S_B$ that define on the space $H$ a structure of a new Hopf algebra that we call the decedent Hopf algebra.  

The author is grateful to Vsevolod Gubarev for his helpful and valuable comments and suggestions.

\section{Basic properties of Rota-Baxter operators on cocommutative Hopf algebras}

Throughout the paper the characteristic of the ground field $F$ is 0. If $A$ is a vector space over $F$ and $\Delta: A\mapsto A$ is a comultiplication on $A$, then we will use the following sumless Sweedler notation for the image of $a\in A$:
$$
\Delta(a)=a_{(1)}\otimes a_{(2)}.
$$

In a Hopf algebra $H=(H,\mu, \Delta, \eta, \epsilon, S)$ we use the following notations:\\
     - $\mu: H\otimes H\mapsto H$ is a multiplication,\\
     - $\Delta: H\mapsto H\otimes H$ is a comultiplication,\\
     - $\eta: F\mapsto H$ is a unit,\\
     - $\epsilon: H\mapsto F$ is a counit,\\
     - $S: H\mapsto H$ is the antipode. 
     
If $(A,\Delta,\epsilon)$ is a coalgebra, then a linear map $\varphi: A\mapsto A$ is called a coalgebra map, if for all $x\in A$:
\begin{gather*}
    \Delta(\varphi(x))=\varphi(x_{(1)})\otimes \varphi(x_{(2)}).\\
    \epsilon(\varphi(x))=\epsilon(x).
\end{gather*}

A Hopf algebra $H$ is called cocommutative if for all $x\in H$ $$x_{(1)}\otimes x_{(2)}=x_{(2)}\otimes x_{(1)}.$$ 
     If $H$ is a cocommutative coalgebra, then the antipode $S:H\mapsto H$ is a coalgebra map. Recall that in arbitrary Hops algebra the antipode $S$ is an algebra antihomomorphism, that is, for all $a,b\in H$ $S(ab)=S(b)S(a)$.
     
     {\bf Definition.} Let $(H,\mu,\eta,\Delta,\epsilon,S)$ be a cocommutative Hopf algebra.  A coalgebra map $B: H\mapsto H$ is called a Rota-Baxter operator on $H$ if for all $x,y\in H$:
\begin{equation}\label{main}B(x)B(y)=B(x_{(1)}B(x_{(2)})yS(B(x_{(3)}))),
\end{equation}
where $\Delta(x)=x_{(1)}\otimes x_{(2)}$. By a Rota-Baxter Hopf algebra we mean a pair $(H,B)$ of a cocommutative Hopf algebra $H$ and a Rota-Baxter operator $B$ on $H$. As an  example on a Rota-Baxter operator on arbitrary cocommutative Hopf algebra one can consider $B=S$, the antipode (see Corollary 2 below).
   
   {\bf Remark.} Note, that if $H$ is a commutative and cocommutative Hopf algebra, then a coalgebra map $B$ is a Rota-Baxter operator if and only if $B$ is an algebra map, that is $B(xy)=B(x)B(y)$ for all $x,y\in H$.

    % In the sequel we will need the following technical
     
{\bf Lemma 1.} Let $H$ be a cocommutative Hops algebra and $B$ be a  Rota-Baxter operator on $H$. Then\\
 (1) If $g\in H$ is a group-like element, then $B(g)$ is also a group-like element.\\
 (2) $B(1)=1$.\\
 (3) If $x\in L$ is a primitive element, then $B(x)$ is also a primitive element. 

{\bf Proof.} 
(1) Let  $g\in H$ be a group-like element. Since $B$ is a coalgebra map, we have
$$
\Delta(B(g))=(B\otimes B)\Delta(g)=B(g)\otimes B(g).
$$
And we have two options: $B(g)=0$ or $B(g)$ is a group-like element of $H$. Since $\epsilon(B(g))=\epsilon(g)=1$, then $B(g)\neq 0$. Therefore, $B(g)$ is a group-like element of $H$.

(2) Since 1 is a group-like element, then so is $B(1)$. Also, by \eqref{main} we have that
$$
B(1)B(1)=B(1B(1)1S(B(1))=B(1).
$$
And since $B(1)$ is inevitable, we get that $B(1)=1$.

(3) Let $x\in H$  be a primitive element. Consider $B(x)$:
$$\Delta(B(x))=(B\otimes B)\Delta(x)=1\otimes B(x)+B(x)\otimes 1.$$
That's mean that $B(x)$ is a primitive element of $H$.

It is well known that if $R:\mathfrak{g}\mapsto \mathfrak{g}$ is a Rota-Baxter operator of weight 1 on a Lie algebra $\mathfrak{g}$, then $-R-id:\mathfrak{g}\mapsto \mathfrak{g}$ is again a Rota-Baxter operator of weight 1 on $\mathfrak{g}$. Similar results for groups was proved in \cite{GLY}: if $G$ is a group and $B$ is a Rota-Baxter operator on $G$, then $\tilde{B}: G\mapsto G$ defined by $\tilde{B}(g)=g^{-1}B(g^{-1})$ is also a Rota-Baxter operator on $G$. For cocommutative Hopf algebras we can generalise these results:

{\bf Proposition 1.} Let $H$ be a cocommutative Hopf algebra and $B$ be a Rota-Baxter operator on $H$. Define $\tilde{B}:H\mapsto H$ as $$\tilde{B}(x)=S(x_{(1)})B(S(x_{(2)})).$$ Then $\tilde{B}$ is also a Rota-Baxter operator on $H$.

{\bf Proof.} Clearly, $\tilde{B}$ is a linear map. Prove that  $\tilde{B}$ is a coalgebra map. Indeed,
\begin{gather*}
\Delta(\tilde{B}(x))=\Delta(S(x_{(1)})B(S(x_{(2)})))=(S(x_{(2)})\otimes S(x_{(1)}))(B(S(x_{(4)}))\otimes B(S(x_{(3)}))=\\
=S(x_{(1)})B(S(x_{(2)}))\otimes S(x_{(3)})B(S(x_{(4)}))=\tilde{B}(x_{(1)})\otimes \tilde{B}(x_{(2)}).
\end{gather*}
In order to prove that $\tilde{B}$ is a Rota-Baxter operator consider
\begin{gather*}
\tilde{B}(x)\tilde{B}(y)=S(x_{(1)})B(S(x_{(2)}))S(y_{(1)})B(S(y_{(2)}))=\\
=S(x_{(1)})\epsilon(x_{(2)})B(S(x_{(3)}))S(y_{(1)})B(S(y_{(2)}))=\\
=S(x_{(1)})B(S(x_{(2)}))S(y_{(1)})\epsilon(B(S(x_{(3)})))B(S(y_{(2)}))=\\
=S(x_{(1)})B(S(x_{(2)}))S(y_{(1)})S(B(S(x_{(3)})))B(S(x_{(4)}))B(S(y_{(2)}))=\\
=hB(S(x_{(3)}))B(S(y_{(2)})),
\end{gather*}
where $h=S(x_{(1)})B(S(x_{(2)}))S(y_{(1)})S(B(S(x_{(3)})))$. For $h$ we have:
\begin{gather*}
h=S(x_{(1)})B(S(x_{(2)}))S(y_{(1)})S(B(S(x_{(3)})))=\tilde{B}(x_{(1)})S(y_{(1)})S(B(S(x_{(2)})))=\\
=\tilde{B}(x_{(1)})S(y_{(1)})S(B(S(x_{(2)})))x_{(3)}S(x_{(4)})=\\
=\tilde{B}(x_{(1)})S(y_{(1)})S(S(x_{(3)})B(S(x_{(2)})))S(x_{(4)})=\tilde{B}(x_{(1)})S(y_{(1)})S(\tilde{B}(x_{(2)}))S(x_{(3)})=\\
=S(x_{(1)}\tilde{B}(x_{(2)})y_{(1)}S(\tilde{B}(x_{(3)}))).
\end{gather*}

Now consider $hB(S(x_{(4)}))B(S(y_{(2)}))$. Using similar arguments as above, we can conclude that:

\begin{gather*}
hB(S(x_{(4)}))B(S(y_{(2)}))=hB(S(x_{(4)})B(S(x_{(5)}))S(y_{(2)})S(B(x_{(6)})))=\\
=hB(S(x_{(4)}\tilde{B}(x_{(5)})y_{(2)}S(\tilde{B}(x_{(6)}))))=\tilde{B}(x_{(1)}\tilde{B}(S(x_{(2)}))y_{(2)}S(\tilde{B}(x_{(3)})).
\end{gather*}
And the proposition is proved.

Another well-known result says that if  a Lie algebra $\mathfrak{g}$ splits into direct sum of two subalgebras $\mathfrak{g}_1$ and $\mathfrak{g}_2$: $\mathfrak{g}=\mathfrak{g}_1\oplus \mathfrak{g}_2$, then the map $R$ defined as $R(x_1+x_2)=-x_2$, where $x_i\in \mathfrak{g}_i$, is a Rota-Baxter operator of weight 1. For groups similar result (\cite{GLY}) says that if a group $G$ can be presented as a product of two subgroups $G_1$ and $G_2$ $G=G_1G_2$ such that  $G_1\cap G_2=\{e\}$, then a map $B$ defined as $B(g_1g_2)=g_2^{-1}$, where $g_i\in G_i$, is a Rota-Baxter operator on $G$. Note, that unlike the Lie algebra case, the inverse of the  projection to the first factor is not a Rota-Baxter operator on $G$.  We can generalise these results for cocommutative Hopf algebras as:

{\bf Proposition 2.} Let $H$ be a cocommutative Hops algebra. Suppose $H_1$ and $H_2$ are two Hopf subalgebras of $H$ and as a Hopf algebra $H=H_1H_2$. Suppose that the product is direct, that is, $H$ is isomorphic to $H_1\otimes _F H_2$ as a vector space. Define a map $B$ as
$$
B(h_1h_2)=\epsilon(h_1)S(h_2),
$$
where $h_i\in H_i$. Then $B$ is a Rota-Baxter operator on $H$.

{\bf Proof.}
 Clearly, $B$ is a well-defined linear map. First we proof that $B$ is a coalgebra map. For $x=\sum h_ig_i$, where $h_i\in H_1$, $g_i\in H_2$ we have:
 \begin{gather*}
 \Delta(B(x))=\sum\limits_i \epsilon(h_i)\Delta(S(g_i))=\sum\limits_i\epsilon(h_i) S(g_{i(2)})\otimes S(g_{i(1)})=\\
 =\sum\limits_i\epsilon(h_{i(1)})S(g_{i(1)})\otimes \epsilon(h_{i(2)})S(g_{i(2)})=(B\otimes B)\Delta(x).
 \end{gather*}
 In order to prove that $B$ satisfies (\ref{main}) consider $x=hg\in H$,  and $y=h'g'$ where $h,h'\in H_1$, $g,g'\in H_2$. We have
 \begin{gather*}
 B(x_{(1)}B(x_{(2)})yS(B(x_{(3)})))=\\
  =B((h_{(1)}g_{(1)})(\epsilon(h_{(2)})S(g_{(2)}))(h'g')(\epsilon(h_{(3)})S(S(g_{(3)})))=\\
 = B(\epsilon(h_{(1)})\epsilon(h_{(2)})h_{(3)}g_{(1)}S(g_{(2)})g'g_{(3)})=B(hh'g'g)=\\
 =\epsilon(hh')S(g'g)=\epsilon(h)\epsilon(h')S(g)S(g')=B(x)B(y).
\end{gather*}
Since $H_1H_2$ is spanned by elements of the form $hg$, the equation (\ref{main}) holds for all $x,y\in H$.

{\bf Remark.} Note that as in the case of groups, the operator $B'$ defined as $B(h_1h_2)=S(h_1)\epsilon(h_2)$ is not a Rota-Baxter operator on $H$ in general. 

{\bf Corollary 1.} If $H$ is a cocommutative Hopf algebra  with the antipode $S$, then $B=S$ is a Rota-Baxter operator on $H$.

\section{Rota-Baxter operators on $F[G]$ and $U(L)$.}

In this section we consider two principle examples of cocommutative Hopf algebra: the group algebra of a group $G$ and the universal enveloping algebra $U(\mathfrak{g})$ of a Lie algebra $\mathfrak{g}$.

{\bf Theorem 1. } Let $(G,B)$ be a Rota-Baxter group. Then  $B$ can be uniquely extended to a  Rota-Baxter  operator $B:F[G]\mapsto F[G]$ on the group algebra $F[G]$. Conversely, if $B$ is a Rota-Baxter operator on $F[G]$, then $B(G)\subset G$ and $(G,B|_G)$ is a Rota-Baxter group, where $B|_G$ is the restriction of $B$ on $G$.

{\bf Proof.}
Since elements of $B$ form a linear basis of $F[G]$,  we can uniquely extend $B$ on $F[G]$ as
$$
B(\sum \alpha_ig_i)=\sum \alpha_i B(g_i).$$

It is easy to see that $B$ is a coalgebra map. We need to check that $(F[G],B)$ is a Rota-Baxter Hopf algebra. Let $x=\sum \alpha_i g_i\in F[G]$, $y\in G$. Then
\begin{gather*}
B(x)B(y)=\sum \alpha_i B(g_i)B(y)= \sum\alpha_i B(g_iB(g_i)yB(g_i)^{-1})\\
=\sum\alpha_iB(g_iB(g_i)hS(B(g_i)))=B(x_{(1)}B(x_{(2)})yS(B(x_{(3)})).
\end{gather*}

And since elements of $G$ form a basis of $F[G]$, the equation (\ref{main}) holds for all $x,y\in F[G]$.

Conversely, let $B$ be a Rota-Baxter operator on the Hopf algebra $F[G]$. By Lemma 1, if $g\in G$, then $B(g)$ is a group-like element. Therefore, $B(g)\in G$ for every $g\in G$.  The rest is obvious.

{\bf Lemma 2.} Let $\mathfrak{g}$ be a Lie algebra and $R$ be a Rota-Baxter operator  of weight 1 on $\mathfrak{g}$. Then the map $R$ can be extended to a linear map $B:U(\mathfrak{g})\mapsto U(\mathfrak{g})$ such that 

1. The restriction of $B$ on $\mathfrak{g}$ is $B|_{\mathfrak{g}}=R$.

2. $B$ satisfies (\ref{main}).

{\bf Proof.} Put $B(1)=1$ and if $x,x_1,\ldots,x_k\in \mathfrak{g}$, $h=x_1x_2\ldots x_k$, then define
\begin{equation}\label{def}
B(xh)=B(x)B(h)-B([B(x),h]).
\end{equation}

First we need to prove that $B$ is well-defined. Consider elements $f,g\in U(\mathfrak{g})$, $x,y\in \mathfrak{g}$. We want to prove that
$B(f(xy-yx-[x,y])g)=0$. 

If $f=1$, then by the definition
\begin{gather*}
B(xyg)=B(x)B(yg)-B([B(x),yg])=B(x)B(yg)-B([B(x),y]g)-B(y[B(x),g])=\\
=B(x)B(y)B(g)-B(x)B([B(y),g])-B([B(x),y])B(g)+B([B([B(x),y]),g])-\\
-B(y)B([B(x),g])+B([B(y),[B(x),g]]).
\end{gather*}
Similarly,
\begin{gather*}
B(xyg)=B(y)B(x)B(g)-B(y)B([B(x),g])-B([B(y),x])B(g)+B([B([B(y),x]),g])-\\
-B(x)B([B(y),g])+B([B(x),[B(y),g]]).
\end{gather*}
And
\begin{gather*}
B(xyg-yxg)=[B(x),B(y)]B(g)-B([B(x),y]+[x,B(y)])B(g)-\\
+B(B([B(x),y]+B(B([x,B(y)]-[[B(x),B(y)],g])=\\
=B([x,y])B(g)-B(B[x,y],g])=B([x,y]g).
\end{gather*}

Suppose that $f=x_1\ldots x_k$ for some $x_i\in \mathfrak{g}$ and use induction on $k$. Denote by $f_1=x_2\ldots x_k$. We have
\begin{gather*}
B(f(xy-yx-[x,y])g)=B(x_1f_1(xy-yx-[x,y])g)=\\
=B(x_1)B(f_1(xy-yx-[x,y])g)-B([B(x_1),f_1(xy-yx-[x,y])g]).
\end{gather*}
By the induction hypothesis, $B(f_1(xy-yx-[x,y])g)=0$. Consider the second summand.
\begin{gather*}
B([B(x_1),f_1(xy-yx-[x,y])g])=
B([B(x_1),f_1](xy-yx-[x,y])g)+\\
+B(f_1[B(x_1),xy-yx-[x,y]]g)+B(f_1(xy-yx-[x,y])[B(x_1),g]).
\end{gather*}
Recall, that $B(x_1)\in \mathfrak{g}$. Then, by the induction hypothesis,
$$B([B(x_1),f_1](xy-yx-[x,y])g)=B(f_1(xy-yx-[x,y])[B(x_1),g])=0.$$ 

Note that 
\begin{gather*}[B(x_1),xy-yx-[x,y]]=
([B(x_1),x]y-y[B(x_1),x]-\\
-[[B(x_1),x],y])+(x[B(x_1),y]-[B(x_1),y]x-[x,[B(x_1),y]])
\end{gather*}
and here we can also use the induction hypothesis to conclude that $B(f_1[B(x_1),xy-yx-[x,y]]g)=0$.

Therefore $B$ is well defined. Now prove that $(U(\mathfrak{g}),B)$ is a Rota-Baxter Hopf algebra. Take $f=x_1\ldots,x_k$, $x_i\in \mathfrak{g}$ and $q\in U(\mathfrak{g})$ and use induction on $k$. If $k=0$ then clearly
$$
B(1)B(q)=B(q)=B(1\cdot B(1)q\cdot S(B(1))).$$

Suppose that $f=xh$ where $h=x_2,\ldots,x_k$. First note that
$$
\Delta^{(2)}(xh)=\sum\limits_{(h)}xh_{(1)}\otimes h_{(2)}\otimes h_{(3)}+h_{(1)}\otimes xh_{(2)}\otimes h_{(3)}+h_{(1)}\otimes h_{(2)}\otimes xh_{(3)}.$$
Then
\begin{gather*}B(f_{(1)}B(f_{(2)})qS(B(f_{(3)})))=B(xh_{(1)}B(h_{(2)})qS(B(h_{(3)})))+B(h_{(1)}B(xh_{(2)})qS(B(h_{(3)})))+\\
+B(h_{(1)}B(h_{(2)})qS(B(xh_{(3)}))).
\end{gather*}

Consider the first term. Using \ref{def} and the induction hypotheses, we have 
\begin{gather*}B(xh_{(1)}B(h_{(2)})qS(B(h_{(3)})))=B(x)B(h_{(1)}B(h_{(2)})qS(B(h_{(3)})))-\\
-B([B(x),h_{(1)}B(h_{(2)})qS(B(h_{(3)}))])=B(x)B(h)B(g)-B([B(x),h_{(1)}B(h_{(2)})qS(B(h_{(3)}))])
\end{gather*}
Similarly,
\begin{gather*}
B(h_{(1)}B(xh_{(2)})qS(B(h_{(3)})))=\\
=B(h_{(1)}B(x)B(h_{(2)})qS(B(h_{(3)})))-B(h_{(1)}B([B(x,h_{(2)}])qS(B(h_{(3)}))))
\end{gather*}
and
\begin{gather*}
B(h_{(1)}B(h_{(2)})qS(B(xh_{(3)})))=\\
=B(h_{(1)}B(h_{(2)})qS(B(x)B(h_{(3)})))-B(h_{(1)}B(h_{(2)})qS(B([B(x),h_{(3)}])))=\\
=-B(h_{(1)}B(h_{(2)})qS(B(h_{(3)}))B(x))-B(h_{(1)}B(h_{(2)})qS(B([B(x),h_{(3)}])))
\end{gather*}

Note that
\begin{gather*}-B([B(x),h_{(1)}B(h_{(2)})qS(B(h_{(3)}))])+B(h_{(1)}B(x)B(h_{(2)})qS(B(h_{(3)})))-\\
-B(h_{(1)}B(h_{(2)})qS(B(h_{(3)}))B(x))=-B([B(x),h_{(1)}]B(h_{(2)})qS(B(h_{(3)}))).
\end{gather*}

Summing up the obtained equations, we get
\begin{gather*}B(f_{(1)}B(f_{(2)})qS(B(f_{(3)})))=\\
=B(x)B(h)B(q)-B([B(x),h_{(1)}]B(h_{(2)})qS(B(h_{(3)})))-\\
-B(h_{(1)}B([B(x,h_{(2)}])qS(B(h_{(3)}))))
-B(h_{(1)}B(h_{(2)})qS(B([B(x),h_{(3)}])))\\
=B(x)B(h)B(q)-B([B(x),h])B(q)=B(xh)B(q).
\end{gather*}
The lemma is proved.

{\bf Lemma 3.} Let $\mathfrak{g}$ be a Lie algebra, $R$ be a Rota-Baxter operator on $\mathfrak{g}$ of weight 1 and $B: U(\mathfrak{g})\mapsto U(\mathfrak{g})$ be the operator from Lemma 2. Then $B$ is a coalgebra map, that is, for all $f\in U(\mathfrak{g})$:
\begin{gather*}
\Delta(B(f))=(B\otimes B)\Delta(f)= B(f_{(1)})\otimes B(f_{(2)}).\\
\epsilon(B(f))=\epsilon(f).
\end{gather*}
{\bf Proof.} First we proof that $B$ preserves the comultiplication.   Take $f=x_1,\ldots x_k$, where $x_i\in \mathfrak{g}$, and use the induction on $k$. The statement is obvious if $f=1$. If $k=1$ we have
$$
\Delta(B(x))=B(x)\otimes 1+1\otimes B(x)=(B\otimes B)(\Delta(x)).$$

Suppose that $f=xh$, where $h=x_2,\ldots,x_k$. Then
$$\Delta(B(xh))=\Delta(B(x)B(h)-B([B(x),h]))=\Delta(B(x))\Delta(B(h))-\Delta(B([B(x),h])$$

By the the induction hypotheses we have:
$$
\Delta(B(h))= B(h_{(1)})\otimes B(h_{(2)})
$$
and since $B(x)\in L$:
$$\Delta(B([B(x),h]))= B([B(x),h_{(1)}])\otimes B(h_{(2)})+B(h_{(1)})\otimes B([B(x),h_{(2)}]).$$
Therefore,
\begin{gather*}
\Delta(B(x))\Delta(B(h))-\Delta(B([B(x),h]))=(B(x)B(h_{(1)}))\otimes B(h_{(2)})+B(h_{(1)})\otimes (B(x)B(h_{(2)}))-\\
-B([B(x),h_{(1)}])\otimes B(h_{(2)})-B(h_{(1)})\otimes B([B(x),h_{(2)}])=\\
=B(xh_{(1)})\otimes B(h_{(2)})+B(h_{(1)})\otimes B(xh_{(2)})=(B\otimes B)(\Delta(xh)).
\end{gather*}

In order to prove that $B$ preserves the counit one can use similar arguments: take $f=x_1,\ldots x_k$ ($x_i\in g$) and use the induction on $k$. The case $k=1$ is trivial. If $f=x_1h$ then
\begin{gather*}
\epsilon(B(xh))=\epsilon(B(x)B(h))-\epsilon(B[B(x),h])=\epsilon(B(x))\epsilon(B(h))-\epsilon(B([B(x),h]))=\\
=\epsilon(x)\epsilon(h)-\epsilon([B(x),h])=\epsilon(xh).\end{gather*}
The lemma is proved.

{\bf Lemma 4.} Let $U(\mathfrak{g})$ be a universal enveloping algebra of a Lie algebra $\mathfrak{g}$, $B$ be a Rota-Baxter operator on $U(\mathfrak{g})$. Then $B(\mathfrak{g})\subset \mathfrak{g}$ and the restriction $R=B|_{\mathfrak{g}}$ is a Rota-Baxter operator of weight 1 on $\mathfrak{g}$.

{\bf Proof.} By Lemma 1, if $x\in \mathfrak{g}$, then $B(x)\in \mathfrak{g}$. Now consider the restriction $R=B|_{\mathfrak{g}}$. For arbitrary $x,y\in \mathfrak{g}$ we have
$$
[R(x),R(y)]=[B(x),B(y)]=B(x)B(y)-B(y)B(x)=B(xy+[B(x),y])-B(yx+[B(y),x])=
$$
$$
=B([B(x),y]+[x,B(y)]+[x,y])=R([R(x),y]+[x,R(y)]+[x,y]).
$$
Therefore, $R$ is a Rota-Baxter operator of weight 1 on the Lie algebra $\mathfrak{g}$.

{\bf Theorem 2.} Rota-Baxter operators of weight 1 on a Lie algebra $\mathfrak{g}$ are in one to one correspondence with Rota-Baxter operators on the universal enveloping algebra $U(\mathfrak{g})$.

{\bf Proof.} Let $R$ be a Rota-Baxter operator  of weight 1 on a Lie algebra $\mathfrak{g}$. It is only left to prove that the extension of $R$ from Lemma 2 is unique.  For this we note, that if $B$ is a Rota-Baxter operator on $U(\mathfrak{g})$, then from \eqref{main} it is follows that
$$
B(x)B(a)=B(xa)+B([B(x),a])$$
for all $x\in \mathfrak{g}$ and $a\in U(\mathfrak{g})$. The rest can be proved using similar arguments  as in Lemma 2.

{\bf Example 1.} Let $\mathfrak{g}=sl_2(F)$, $x,h,y$ be a basis of $sl_2(F)$ with the following table of multiplication:
$$
[h,x]=2x,\quad [h,y]=-2y,\quad [x,y]=h.
$$
Consider a map $R$ defined as
$$
R(x)=0,\quad R(h)=-\frac{h}{2},\quad R(y)=-y.
$$
Then $R$ is a Rota-Baxter operator of weight 1 on $\mathfrak{g}$ (\cite{GME}). Consider the extension $B$ of $R$ on $U(\mathfrak{g})$. Let $a=xb$, where $b\in U(\mathfrak{g})$. By the definition of $B$:
$$
B(a)=B(xb)=B(x)B(h)-B([B(x),b])=0.
$$
Similarly, if $a=yb$, we have:
$$
B(a)=B(yb)=B(y)B(b)-B([B(y),b])=-yB(b)+B([y,b])=-yB(b)+B(yb)-B(by).
$$
Therefore, $B(by)=-yB(b)$. Monomials  $x^ih^jy^k$ form a linear basis of $U(\mathfrak{g})$. For them we proved that 
$$
B(x^ih^jy^k)=\left \{ \begin{matrix} 0,\ \text{if}\  i>0\\
(-1)^{j+k}\frac{y^kh^j}{2^j},\ \ \text{if}\ i=0.
\end{matrix}\right.
$$
\section{ The descendent Hopf algebra}

{\bf Definition \cite{Val}.} Let $(\mathfrak{g},[,])$ be a Lie algebra and $\cdot$ is a bilinear operation on L. If for all $x,y,z\in g$:
\begin{gather*}
[x,y]\cdot z=(y\cdot x)\cdot z-y\cdot(x\cdot z)-(x\cdot y)\cdot z+x\cdot(y\cdot z),\\
x\cdot[y,z]=[x\cdot y,z]+[y,x\cdot z],
\end{gather*}
then $(\mathfrak{g},[,],\cdot)$ is called a Post-Lie algebra.

Given a Post-Lie algebra $(\mathfrak{g},[,],\cdot)$, one can define new Lie bracket on $\mathfrak{g}$ by the formula
\begin{equation}\label{new}
\{x,y\}=x\cdot y-y\cdot x+[x,y].
\end{equation}

If $\mathfrak{g}$ is a Lie algebra and $R$ is a Rota-Baxter operator of weight 1 on $\mathfrak{g}$, then one can define the following structure of a post-Lie algebra:
$$
x\cdot y=[R(x),y]$$
for all $x,y\in \mathfrak{g}$ \cite{Bai}. In this case the multiplication \eqref{new} is equal to
$$
\{x,y\}=[R(x),y]+[x,R(y)]+[x,y].
$$

{\bf Definition \cite{GLY}.} If $R$ is a Rota-Baxter operator of weight 1 on a Lie algebra $\mathfrak{g}$, then the pair $(\mathfrak{g},\{,\})$ is called the descendent Lie algebra of the Rota-Baxter Lie algebra $(\mathfrak{g},R)$.  

In the same paper it was proved that

{\bf Statement 1 \cite{GLY}.} If $(G,B)$ is a Rota-Baxter group, then $G$ with the product 
$$
g*h=gB(g)hB(g)^{-1}
$$
is again a group called the descendent group $G_B$ of the Rota-Baxter group $(G,B)$. The inverse of $g\in G_B$ in the descendent group is equal to $B(g)^{-1}g^{-1}B(g)$.

Let $(\mathfrak{g},[,],\cdot)$ be a post-Lie algebra. In \cite{FLMK} it was proved, that there is a unique extension of the post-Lie product $\cdot$ on the universal enveloping algebra $U(\mathfrak{g})$ given by
\begin{gather}
1\cdot f=f, \label{e1} \\
xf\cdot g= x\cdot(f\cdot g)-(x\cdot f)\cdot g, \label{e2} \\
f\cdot(gh)= (f_{(1)}\cdot g)(f_{(2)}\cdot h) \label{e3}
\end{gather}
for all $x\in \mathfrak{g}$, $f,g\in U(\mathfrak{g})$.

Define a new multiplication on $U(\mathfrak{g})$ by
$$f*g= f_{(1)}(f_{(2)}\cdot g)$$
where $f,g\in U(\mathfrak{g})$. In the same paper it was proved that $(U(\mathfrak{g}),*)$ is isomorphic to the universal enveloping algebra of $(\mathfrak{g},\{,\})$. 

{\bf Proposition 3.} Let $(\mathfrak{g},[,])$ be a Lie algebra, $R$ be a Rota-Baxter operator on $\mathfrak{g}$ of weight 1, $(\mathfrak{g},[,],\cdot)$ be the corresponding post-Lie algebra and $(U(\mathfrak{g}),B)$ be the enveloping Rota-Baxter algebra of $(\mathfrak{g},R)$. Then the extension $\cdot$ of the post-Lie product on $U(\mathfrak{g})$ can be defined as
$$f\cdot g= B(f_{(1)})gS(B(f_{(2)})).
$$
{\bf Proof.} Since the extension defined by (\ref{e1})-(\ref{e3}) is unique, it is enough   to prove that our product satisfies  (\ref{e1})-(\ref{e3}).

Let $x\in \mathfrak{g}$, $f,g,h\in U(\mathfrak{g})$. The first equation is obvious. Consider (\ref{e2}). We have
\begin{gather*}
x\cdot(f\cdot g)-(x\cdot f)\cdot g= x\cdot (B(f_{(1)})gS(B(f_{(2)})))-[B(x),f]\cdot g=\\
=[B(x),B(f_{(1)})gS(B(f_{(2)})]-B([B(x),f_{(1)}])gS(B(f_{(2)}))-\\
-B(f_{(1)})gS([B(x),B(f_{(2)})])
=B(xf_{(1)}gS(B(f_{(2)}))+B(f_{(1)})gS(B(xf_{(2)}))=(xf)\cdot g.
\end{gather*}

Consider $(\ref{e3})$. We have
\begin{gather*} (f_{(1)}\cdot g)(f_{(2)}\cdot h)= (B(f))_{(1)}gS((B(f))_{(2)})(B(f))_{(3)}hS((B(f))_{(4)})=\\
= (B(f))_{(1)}g\epsilon ((B(f))_{(2)}hS((B(f))_{(3)}))= (B(f))_{(1)}ghS((B(f))_{(2)})=f\cdot (gh).\\
\end{gather*}
The proposition is proved.

{\bf Corollary 2.} Let $(\mathfrak{g},[,])$ be a Lie algebra, $R$ be a Rota-Baxter operator of weight 1 on $\mathfrak{g}$ and $B$ be the extension of $R$ on $U(\mathfrak{g})$ from Lemma 2. Define new multiplication $*: U(\mathfrak{g})\otimes U(\mathfrak{g})\mapsto U(\mathfrak{g})$ as
$$
f*g= f_{(1)}B(f_{(2)})gS(B(f_{(3)})).
$$
Then $(U(\mathfrak{g}),*,\Delta,\eta,\epsilon)$ is a bialgebra  isomorphic to the universal enveloping algebra of the Lie algebra $(g,\{,\})$, where product $\{,\}$ is defined as
$$
\{x,y\}=[R(x),y]+[x,R(y)]+[x,y]
$$
for all $x,y\in \mathfrak{g}$.

{\bf Remark.} As we will see in Theorem 4, the antipode $S_B$ on $(U(\mathfrak{g}),*,\Delta,\eta,\epsilon)$ is  defined as
$$S_B(x)=S(B(x_{(1)}))S(x_{(2)})B(x_{(3)})
$$
for all $x\in U(\mathfrak{g})$. Here, $S$ is the "old" antipode of the universal enveloping algebra  $U(\mathfrak{g})$ of $\mathfrak{g}$.

We want to generalise Corollary 2 to arbitrary cocommuative Rota-Baxter Hopf algebra.  Let $(H,\mu, \Delta, \eta, \epsilon, S)$ be a cocommutative Hopf algebra and $B$ a Rota-Baxter operator on $H$. Define new operation on $H$: for all $x,y\in H$ put
$$x*y=x_{(1)}B(x_{2})yS(B(x_{3})).
$$

{\bf Proposition 4}. We have that
\begin{gather*}
\Delta(x*y)=(x_{(1)}*y_{(1)})\otimes (x_{(2)}*y_{(2)}),\\
\epsilon(x*y)=\epsilon(x)\epsilon(y).
\end{gather*}

{\bf Proof.}
Consider $x*y=x_{(1)}B(x_{2})yS(B(x_{3}))$. Since $\Delta$ preserves multiplication, we have
$$
\Delta(x*y)=\Delta(x_{(1)}B(x_{2})yS(B(x_{3})))=
$$
$$
=(x_{(1)}\otimes x_{(2)})(B(x_{(3)})\otimes B(x_{(4)}))(y_{(1)}\otimes y_{(2)})(S(B(x_{(6)}))\otimes S(B(x_{(5)})))
$$
Since the comultiplication is cocommutative, we can rewrite the last term as
$$
x_{(1)}B(x_{(2)})y_{(1)}S(B(x_{(3)}))\otimes x_{(4)}B(x_{(5)})y_{(2)}S(B(x_{(6)}))=x_{(1)}*y_{(1)}\otimes x_{(2)}*y_{(2)}.
$$

Now consider
\begin{gather*}
\epsilon(x*y)=\epsilon(x_{(1)}B(x_{2})yS(B(x_{3})))=\epsilon(x_{(1)})\epsilon(B(x_{2}))\epsilon(y)\epsilon(S(B(x_{3})))=\\
=\epsilon(x_{(1)})\epsilon(B(x_{2})S(B(x_{3})))\epsilon(y)=
\epsilon(x_{(1)})\epsilon(B(x_{(2)})\epsilon(y)=\epsilon(x_{(1)}\epsilon(x_{(2)}))\epsilon(y)=\\
=\epsilon(x)\epsilon(y).
\end{gather*}

Define a linear map $S_B: H\mapsto H$ as
$$
S_B(x)=S(B(x_{(1)}))S(x_{(2)})B(x_{(3)})$$
for all $x\in H$. We will need the following 

{\bf Proposition 5.} For all $x\in H$ we have
$$
\epsilon(x)1=B(x_{(1)})B(S_B(x_{(2)})).
$$

{\bf Proof.}
Indeed, 
\begin{gather*}
B(x_{(1)})B(S_B(x_{(2)}))=B(x_{(1)}B(x_{(2)})S_B(x_{(3)})S(B(x_{(4)})))=\\
=B(x_{(1)}[B(x_{(2)})S(B(x_{(3)}))]S(x_{(4)})[B(x_{(5)})S(B(x_{(6)}))])=\\
B(x_{(1)}\epsilon (B(x_{(2)}))S(x_{(3)})\epsilon(B(x_{(4)})))=
B(x_{(1)}S(x_{(2)}))=\epsilon(x)B(1)=\epsilon(x)1.
\end{gather*}

{\bf Theorem 4.}
$H_B=(H,*,\Delta, \eta, \epsilon, S_B)$ is a cocommutative Hopf algebra.

{\bf Proof.} Note that since $B$ is a Rota-Baxter operator on $H$, we have that
$$B(x*y)=B(x)B(y)
$$
for all $x,y\in H$. 

First we proof that $(H,*)$ is an associative algebra. Indeed, take $x,y,z\in H$. Using Proposition 4 and Proposition 5, we have
\begin{gather*}
(x*y)*z=(x_{(1)}*y_{(1)})B(x_{(2)}*y_{(2)})zS(B(x_{(3)}*y_{(3)})=\\
=x_{(1)}B(x_{(2)})y_{(1)}[S(B(x_{(3)}))B(x_{(4)})]B(y_{(2)})zS(B(y_{(3)}))S(B(x_{(5)}))=\\
=x_{(1)}B(x_{(2)})y_{(1)}[\epsilon(B(x_{(3)}))]B(y_{(2)})zS(B(y_{(3)}))S(B(x_{(4)}))=\\
=x_{(1)}B(x_{(2)})y_{(1)}B(y_{(2)})zS(B(y_{(3)}))S(B(x_{(3)}))=
x*(y*z).
\end{gather*}

Since $B(1)=1$, it is easy to see that $1*x=x*1=x$ for all $x\in H$.

By Proposition 4, $(H,*,\eta, \Delta,\epsilon)$ is a bialgebra. It is left to proof that $S_B$ is the antipode of  $(H,*,\Delta, \eta, \epsilon)$. We need to prove that for all $x\in H$:
$$
x_{(1)}*S_B(x_{(2)})=S_B(x_{(1)})*x_{(2)}=\epsilon(x)1.
$$

Direct computation shows
\begin{gather*}
x_{(1)}*S_B(x_{(2)})=x_{(1)}[B(x_{(2)})S(B(x_{(3)}))]S(x_{(4)})[B(x_{(5)})S(B(x_{(6)}))]=\\
=x_{(1)}\epsilon(x_{(2)})S(x_{(3)})\epsilon(x_{(4)})=x_{(1)}S(x_{(2)})=\epsilon(x)1.
\end{gather*}

Note that $B(x_{(1)}*S_B(x_{(2)}))=B(x_{(1)})B(S_B(x_{(2)}))$. Then we get the equality
\begin{equation}\label{anteq}
B(x_{(1)})B(S_B(x_{(2)}))=\epsilon(x)1.
\end{equation}

Also, we have that
\begin{gather*}
B(S_B(x_{(1)}))B(x_{(2)})=B(S_B(x_{(1)}))\epsilon(B(x_{(2)}))B(x_{(3)})=\\
=\epsilon(B(x_{(1)}))B(S_B(x_{(1)}))B(x_{(3)})=S(B(x_{(1)}))[B(x_{(2)})B(S_B(x_{(3)}))]B(x_{(4)})=\\
= S(B(x_{(1)}))\epsilon(x_{(2)})B(x_{(3)})=\epsilon(x)1.
\end{gather*}
That is, we proved that
\begin{equation}\label{anteq2}
   B(S_B(x_{(1)}))B(x_{(2)})=\epsilon(x)1. 
\end{equation}

Now consider the second equality. Using (\ref{anteq}) and (\ref{anteq2}) we compute:  
\begin{gather*}
S_B(x_{(1)})*x_{(2)}=S_B(x_{(1)})B(S_B(x_{(2)}))x_{(3)}S(B(S_B(x_{(4))}))=\\
=S(B(x_{(1)}))S(x_{(2)})[B(x_{(3)})B(S_B(x_{(4)}))]x_{(5)}S(B(S_B(x_{(6))}))=\\
=S(B(x_{(1)}))S(x_{(2)})[\epsilon(x_{(3)})]x_{(4)}S(B(S_B(x_{(5))}))=\\
=S(B(x_{(1)}))S(x_{(2)})x_{(3)}S(B(S_B(x_{(4))}))=S(B(x_{(1)}))S(B(S_B(x_{(2))}))=\\
=S(B(S_B(x_{(2)})B(x_{(1)}))=\epsilon(x)1.
\end{gather*}
And the theorem is proved.

{\bf Proposition 6.} $B$ is a homomorphism of Hopf algebras $H$ and $H_B$ and is a Rota-Baxter operator  on the Hopf algebra $H_B$.

{\bf Proof.} First we note, that by the definition of the product $*$ and by (\ref{main}), for all $x,y\in H$ we have that: $B(x*y)=B(x)B(y)$. It is left to to proof that $B\circ S_B=S\circ B$. For this note that 
$$
S(B(x_{(1)}))B(x_{(2)})=B(x_{(1)})S(B(x_{(2)}))=\epsilon(B(x))=\epsilon(x)1.
$$
This means that the map $S\circ B$ is the inverse for the map $B$ in $(End(B),\star,\eta\circ\epsilon)$ where $\star:End(B)\mapsto End(B)$ is the convolution product defined as
$$
(f\star g)(x)=f(x_{(1)})g(x_{(2)})$$
for all $f,g\in \End(H),\ x\in H$.

On the other hand, for every $x\in H$:
$$
B(S_B(x_{(1)}))B(x_{(2)})=B(S_B(x_{(1)})*x_{(2)})=B(\epsilon(x)1)=\epsilon(x)1.
$$
Similarly, $B(x_{(1)})B(S_B(x_{(2)}))=\epsilon(x)1$ and $B\circ S_B$ is also the inverse for $B$ in $(\End(B),\star,\eta\circ\epsilon)$. Therefore, $B\circ S_B=S\circ B$ and $B$ is a homomorphism of Hopf algebras $H$ and $H_B$.

For the second statement we compute:
\begin{gather*}
B(x_{(1)}*B(x_{(2)})*y*S(B(x_{(3)})))=B(x_{(1)})B(B(x_{(2)}))B(y)B(S(B(x_{(3)})))=\\
=B(x)*B(y).
\end{gather*}

By analogy with Lie algebras and groups, we may give the following 

{\bf Definition.} The Hopf algebra $H_B$ is called the descendent Hopf algebra of the Rota-Baxter Hopf algebra $(H,B)$.

{\bf Remark.} Corollary 2 says that any descendent Hopf algebra of the universal enveloping algebra $U(\mathfrak{g})$ of a Lie algebra $\mathfrak{g}$ is the universal enveloping algebra of the corresponding descendent Lie algebra.  And it is easy to see that if $H=F[G]$ is the group algebra of a group $G$, then any descendent Hopf algebra $H_B$ is the group algebra of the correspondent descendent group $G_B$.

\section*{Acknowledgements}

The work was supported by Russian Scientific Fond (project N 19-11-00039).

\noindent Maxim Goncharov \\
Novosibirsk State University \\
Sobolev Institute of Mathematics \\
Novosibirsk, Russia \\
e-mail: goncharov.gme@gmail.com

\end{document}